\documentclass[12pt,a4paper]{article}
\usepackage{amssymb,amsfonts,amsthm,amsmath}

\usepackage{epsf}
\usepackage{graphicx}
\usepackage{placeins}
\usepackage{setspace}

\def\hang{\hangindent\parindent}
 \def\rf{\par\noindent\hang}

\begin{document}
\baselineskip=20pt

\begin{center}
\Large{{\bf A comparison of Bayesian and frequentist interval estimators
in regression that utilize uncertain prior information}}
\end{center}

\smallskip

\begin{center}
\large{
{\bf Paul Kabaila$^*$ and Gayan Dharmarathne}}
\end{center}

\begin{center}
{\sl Department of Mathematics and Statistics,
La Trobe University, Victoria 3086, Australia}
\end{center}


\noindent {\bf ABSTRACT}

\medskip

\noindent Consider a linear regression model with regression parameter $\boldsymbol{\beta}$ and normally distributed errors.
Suppose that the parameter of interest is $\theta = \boldsymbol{a}^T \boldsymbol{\beta}$
where $\boldsymbol{a}$ is a specified vector.
Define the parameter $\tau = \boldsymbol{c}^T \boldsymbol{\beta} - t$ where  $\boldsymbol{c}$ and $t$
are specified and $\boldsymbol{a}$ and $\boldsymbol{c}$ are linearly independent.
Also suppose that we have uncertain prior information that $\tau = 0$.
Kabaila and Giri, 2009, JSPI, describe a new frequentist $1-\alpha$ confidence interval for
$\theta$ that utilizes this uncertain prior information.
We compare this confidence interval with Bayesian $1-\alpha$ equi-tailed and shortest credible intervals
for $\theta$ that result from a prior density for $\tau$ that is a mixture of a
rectangular ``slab'' and a Dirac delta function ``spike'', combined
with noninformative prior densities for the other parameters of the model. We show that these frequentist and Bayesian
interval
estimators depend on the data in very different ways. We also consider some close variants of this prior distribution
that lead to Bayesian and frequentist interval estimators with greater similarity. Nonetheless, as we show,
substantial differences between these interval estimators remain.

\bigskip


\vspace{0.7cm}

\noindent $^*$Corresponding author. Tel.: +61 3 9479 2594, fax: +61 3 9479 2466.

\noindent {\sl E-mail address:} P.Kabaila@latrobe.edu.au (Paul Kabaila).

\newpage

\noindent {\large{\bf 1. Introduction}}

\medskip

Consider the linear regression model $\boldsymbol{Y} = \boldsymbol{X} \boldsymbol{\beta} + \boldsymbol{\varepsilon}$,
where $\boldsymbol{Y}$ is a random $n$-vector of responses, $\boldsymbol{X}$ is a known $n \times p$ matrix with linearly
independent columns, $\boldsymbol{\beta} = (\beta_1,\ldots, \beta_p)$ is an unknown parameter vector ($p \ge 2$) and
$\varepsilon \sim N(0, \sigma^2 I_n)$ where $\sigma^2$ is an unknown positive parameter.
Suppose that the parameter of interest is $\theta = \boldsymbol{a}^T \boldsymbol{\beta}$
where $\boldsymbol{a}$ is a specified
$p$-vector ($\boldsymbol{a} \ne \boldsymbol{0}$).
The inference of interest is an interval estimator for $\theta$.
Define the parameter $\tau = \boldsymbol{c}^T \boldsymbol{\beta} - t$ where the vector $\boldsymbol{c}$ and the number $t$
are specified and $\boldsymbol{a}$ and $\boldsymbol{c}$ are linearly independent.
Also suppose that previous experience with similar data sets and/or
expert opinion and scientific background suggest that $\tau = 0$.
In other words, suppose that we have uncertain prior information that $\tau = 0$.
Kabaila and Giri (2013) describe six examples of this scenario. These include a $2^k$ factorial experiment with two or more replicates,
where the parameter of interest $\theta$ is a specified contrast and the uncertain prior information is that the highest order
interaction is zero.
For clarity of comparison of the interval estimators considered,
we assume that Var$(\hat{\theta}) = \sigma^2$ and
Var$(\hat{\tau}) = \sigma^2$. In Appendix A, we show that this can be achieved by appropriate scaling
and that there is no loss of generality, as far as the purposes of the paper are concerned.

The uncertain prior information about $\tau$ can be utilized in the construction of the interval estimator
for $\theta$ in two ways: Bayesian and frequentist.
A Bayesian $1-\alpha$ credible interval for $\theta$ that utilizes the uncertain prior
that $\tau = 0$ is obtained by using an informative prior for $\tau$, combined with
noninformative priors for the other parameters in the model.
A frequentist confidence interval for $\theta$ is said to be a $1-\alpha$ confidence
interval if it has infimum coverage probability $1-\alpha$.
We assess a $1-\alpha$ confidence
interval $J$ by its scaled expected length, defined to be the ratio (expected length of $J$)/(expected
length of the standard $1-\alpha$ confidence interval for $\theta$).
Farchione and Kabaila (2008), Kabaila (2009) and Kabaila and Giri (2009, 2013) define a
 frequentist $1-\alpha$ confidence interval for $\theta$ to be one that
 utilizes the uncertain prior information
that $\tau = 0$ if it has the following properties:
(a) the scaled expected length of this interval is substantially
less than 1 when $\tau = 0$, (b) the maximum (over the parameter space)
of the scaled expected length is not too large and (c) this
confidence interval reverts to the standard $1-\alpha$ confidence interval
when the data happen to strongly contradict the prior information.
The strong admissibility of the standard $1-\alpha$ confidence interval
(Kabaila, Giri and Leeb, 2010) implies that the maximum value of the scaled expected
length must be greater than 1.

Kabaila and Giri (2009, 2013) describe a frequentist
$1-\alpha$ confidence interval for $\theta$ that utilizes the uncertain prior information
that $\tau=0$.
For brevity, we refer to this as the KG $1-\alpha$ confidence interval.
It is important to compare the KG $1-\alpha$ confidence interval with
Bayesian $1-\alpha$ credible intervals for
$\theta$ that utilize the uncertain prior information
that $\tau = 0$.
Our assumption is that we have no prior information about $\boldsymbol{\beta}$, apart from the
uncertain prior information that $\tau=0$. Because this prior information is so precisely targeted, it is
inappropriate to use a $g$-prior (Zellner, 1986) for the construction of a Bayesian credible interval
for $\theta$. Even so, there is a multitude of possible informative prior distributions
for $\tau$, each leading to a different Bayesian $1-\alpha$ credible interval for $\theta$.
In the present paper, we deal exclusively with the following prior distribution and its close variants.
This prior results from an improper prior density for $\tau$ that
consists of a mixture of an infinite rectangular unit-height ``slab'' and a Dirac delta  function ``spike'', combined
with noninformative prior densities for the other parameters of the model.
This prior belongs to the class of `slab and spike' (or `spike and slab') priors that are a mixture of a Dirac delta function
``spike'' and a density function that is symmetric about 0 and achieves its maximum at 0. This class of priors
is widely used for Bayesian variable selection, see e.g.
Mitchell and Beauchamp (1988), Chipman, George and McCulloch (2001), Section 7.2 of Miller (2002)
and O'Hara and Sillanp\"a\"a (2009). This class of priors is also used for estimation under the assumption of possible sparsity,
see e.g. Johnstone and Silverman (2004, 2005). Variable selection may be an end in itself, e.g. in genomic studies that aim
to predict disease outcome. However, in scenarios such as those considered in the present paper, any variable selection is just a
preliminary step to finding an interval estimator for $\theta$. For Bayesian interval estimation, it makes sense to use the same prior
that one would use for Bayesian variable selection.

So, in Section 3, we suppose that the prior density is
$\pi(\theta, \tau, \sigma^2) = \big(\xi \delta(\tau) + (1-\xi) \big) \sigma^{-2}$,
where $\delta$ denotes the Dirac delta function and $\xi \in [0,1]$. Although this is an improper prior density, the marginal posterior
distribution of $\theta$ is a well-behaved proper distribution. The parameter $\xi$ specifies the strength of
the prior belief that $\tau=0$. The strength of this prior belief increases with
increasing $\xi$, with $\xi=0$ corresponding to no prior information about $\tau$ and $\xi=1$ corresponding
to certainty that $\tau$ is 0.
An attractive feature of this prior density is that the Bayesian  $1-\alpha$ highest posterior density (HPD) and
equi-tailed credible intervals for $\theta$ are identical to the usual frequentist
$1-\alpha$ confidence interval in the two
extreme cases that (a) $\tau$ is known to be 0 (i.e. $\xi = 1$) and (b) there is no prior information about $\tau$
(i.e. $\xi = 0$).

In Section 3, we show that, for $\xi \in (0,1)$, the
$1-\alpha$ HPD credible set for $\theta$ may consist of a union of two disjoint intervals.
This is because, as illustrated by Figure 1, the marginal posterior density of $\theta$
may be bimodal.
We therefore focus on Bayesian $1-\alpha$ equi-tailed and
shortest credible intervals for $\theta$.

Let $\hat{\boldsymbol{\beta}}$ denote the least squares estimator of $\boldsymbol{\beta}$. Now let
$\hat{\theta} = \boldsymbol{a}^T \hat{\boldsymbol{\beta}}$,
$\hat{\tau} = \boldsymbol{c}^T \hat{\boldsymbol{\beta}} - t$ and
$\hat \sigma^2 = (\boldsymbol{Y} - \boldsymbol{X} \hat{\boldsymbol{\beta}})^T (\boldsymbol{Y} - \boldsymbol{X} \hat{\boldsymbol{\beta}})/(n-p)$.
We describe both frequentist and Bayesian interval estimators for $\theta$ using the {\sl scaled half-length}, defined to be
\begin{equation*}
\frac{\text{length of interval}}{2 \, \hat{\sigma}},
\end{equation*}
and the {\sl scaled offset}, defined to be
\begin{equation*}
\frac{(\text{centre of interval}) - \hat{\theta}}{\hat{\sigma}}.
\end{equation*}

For the KG $1-\alpha$ confidence interval,
both the scaled half-length and the scaled offset are functions of $\hat{\tau} / \hat{\sigma}$. This
makes sense
because $|\hat{\tau}| / \hat{\sigma}$ is a frequentist measure of
the extent to which that data are inconsistent with the uncertain prior information that $\tau=0$.
We show in Section 3 (where we suppose that the prior density is
$\pi(\theta, \tau, \sigma^2) = (\xi \delta(\tau) + (1-\xi) ) \sigma^{-2}$) that, in sharp contrast to this, for Bayesian
$1-\alpha$ equi-tailed and shortest credible intervals both the scaled half-length and
the scaled offset are functions of $(\hat{\sigma}, \hat{\tau} / \hat{\sigma})$, for all $\xi \in (0,1)$.
This is illustrated in Figures 2 and 3. Figure 2 shows graphs of the scaled offset and the scaled half-length
of a Bayesian 0.95 equi-tailed credible interval, as functions of $\hat{\tau} / \hat{\sigma}$,
for $\hat{\sigma}=1$ (solid line) and $\hat{\sigma}=10$ (dashed line). Figure 3 shows graphs of the scaled offset and the scaled half-length
of the Bayesian 0.95 shortest credible interval, as functions of $\hat{\tau} / \hat{\sigma}$,
for $\hat{\sigma}=1$ (solid line) and $\hat{\sigma}=10$ (dashed line).
In other words, for the prior distribution considered in Section 3, we show that
the KG $1-\alpha$ confidence interval depends on the data
in a very different way from
the Bayesian $1-\alpha$ equi-tailed and shortest credible intervals.

In Section 4, we consider some close variants of the informative prior distribution considered in Section 3,
that lead to Bayesian and KG frequentist interval estimators with greater similarity,
in that both the scaled half-length and the scaled offset are functions of $\hat{\tau} / \hat{\sigma}$.
This is illustrated in Figures 4 and 5.
Nonetheless, there are still very substantial differences
between the Bayesian $1-\alpha$ equi-tailed and shortest credible intervals and the KG frequentist $1-\alpha$ confidence interval.
Our conclusion is that the KG $1-\alpha$ confidence interval and Bayesian $1-\alpha$ equi-tailed and shortest credible intervals
for $\theta$ that utilize the uncertain prior information that $\tau = 0$ are
different for the informative prior distributions considered in both Sections 3 and 4.

\bigskip

\noindent {\large{\bf 2. Brief description of the KG $\boldsymbol{1-\alpha}$ confidence interval}}

\medskip

The standard $1-\alpha$ confidence interval for $\theta$ is
$I = \big [ \hat \theta - t(m)  \hat
\sigma, \quad \hat \theta + t(m)  \hat
\sigma \big ]$,
where $m=n-p$ and the quantile $t(q)$ is defined by $P(-t(q) \le T \le t(q)) = 1-\alpha$ for $T \sim t_q$.
The following is a brief description of the KG $1-\alpha$ confidence interval.
Suppose that $b: \mathbb{R} \rightarrow \mathbb{R}$ is a continuous odd function and
$s: \mathbb{R} \rightarrow (0, \infty)$ is a continuous even function. Also suppose that
$b(x)=0$ for all $|x| \ge d$ and $s(x)=t(m)$ for all $x \ge d$,
where $d$ is a (sufficiently large) specified positive number.
For each $b$ and $s$,
define the following confidence interval for $\theta$:
\begin{align}
\label{J(b,s)}
J(b,  s)
= \left [ \hat{\theta} - \hat{\sigma} \, b \left(\frac{\hat{\tau}}{\hat{\sigma}} \right)
- \hat{\sigma} \, s \left(\frac{\hat{\tau}}{\hat{\sigma}} \right), \,
\hat{\theta} - \hat{\sigma} \, b \left(\frac{\hat{\tau}}{\hat{\sigma}} \right)
+ \hat{\sigma} \, s \left(\frac{\hat{\tau}}{\hat{\sigma}} \right)
 \right ].
\end{align}
For this interval estimator, the scaled half-length is $s(\hat{\tau}/\hat{\sigma})$
and the scaled offset is $-b(\hat{\tau}/\hat{\sigma})$. In other words, both the
scaled half-length and the scaled offset are functions of $\hat{\tau}/\hat{\sigma}$.
The statistic $|\hat{\tau}|/\hat{\sigma}$ is the usual frequentist test
statistic for testing the null hypothesis $H_0: \tau=0$ against the alternative hypothesis $H_1: \tau \ne 0$.
This implies that the confidence interval $J(b,  s)$ reverts to the standard $1-\alpha$ confidence interval
$I$ when the data happen to strongly contradict the uncertain prior information that $\tau = 0$.
Kabaila and Giri (2009) and Kabaila and Giri (2013) describe two methods for the computation of smooth functions $b$ and $s$
such that $J(b, s)$ is a $1-\alpha$ confidence interval for $\theta$ that utilizes the uncertain prior
information that $\tau = 0$, in the sense described in the second paragraph of the introduction.
This computation is carried out
by the statistician {\sl prior} to looking at the observed
response vector $\boldsymbol{y}$.

\bigskip

\noindent {\large{\bf 3. Comparison of the frequentist and Bayesian interval estimators for
the prior density $\boldsymbol{\pi(\theta, \tau, \sigma^2)  = \big(\xi \delta(\tau) + (1-\xi) \big) {\sigma}^{-2}}$}}

\medskip

In this section, our aim is to compare this KG $1-\alpha$ confidence interval with Bayesian
$1-\alpha$ credible intervals for $\theta$ that result from the improper prior density
$\xi \delta(\tau) + (1-\xi)$
for $\tau$, combined with noninformative prior
distributions for the other parameters of the model.
In Appendix B, we define the parameter vector $\boldsymbol{\chi}$,
which has dimension $p-2$. In this appendix, we use
sufficiency to reduce the data to
$\big(\hat{\theta}, \hat{\tau}, \hat{\boldsymbol{\chi}}, \hat{\sigma}^2\big)$,
where $\hat{\boldsymbol{\chi}}$ is the least squares estimator of $\boldsymbol{\chi}$.
Under the sampling model,
the random vectors $\big(\hat{\theta}, \hat{\tau} \big)$,
$\hat{\boldsymbol{\chi}}$ and $\hat{\sigma}^2$ are independent,
\begin{equation}
\label{model_thetahat_tauhat}
\left[\begin{matrix} \hat{\theta}\\ \hat{\tau} \end{matrix}
\right] \sim N \left ( \left[\begin{matrix} \theta \\ \tau \end{matrix}
\right], \sigma^2 \left[\begin{matrix} 1 \quad \rho\\ \rho \quad 1 \end{matrix}
\right] \right ),
\end{equation}
$\hat{\boldsymbol{\chi}} \sim N \big(\boldsymbol{\chi}, \sigma^2 \boldsymbol{I}_{p-2} \big)$  and
$m \hat{\sigma}^2/ \sigma^2 \sim \chi^2_m$. Throughout the Bayesian analysis in this paper,
we suppose that the prior distributions of
$(\theta, \tau, \sigma^2)$ and $\boldsymbol{\chi}$ are independent and that
the components of $\boldsymbol{\chi}$ have independent uniform prior distributions. As shown in Appendix B,
the marginal posterior distribution of $(\theta, \tau, \sigma^2)$ is the same as the posterior
distribution of $(\theta, \tau, \sigma^2)$ based on the reduced data
$\big(\hat{\theta}, \hat{\tau}, \hat \sigma^2 \big)$ and the sampling model that
$\big(\hat{\theta}, \hat{\tau}\big)$ and $\hat \sigma^2$ are independent random vectors, $\big(\hat{\theta}, \hat{\tau}\big)$
has the distribution \eqref{model_thetahat_tauhat} and $m \hat{\sigma}^2/ \sigma^2 \sim \chi^2_m$.
It is this reduced data and the corresponding sampling model that we use from now on for our Bayesian analysis.
In addition,
throughout our Bayesian analysis,
we suppose that the prior distributions
of $\theta$ and $(\tau,\sigma^2)$ are independent, that $\theta$ has a uniform prior distribution over the real line
and $\sigma^2$ has the improper prior density $1/ \sigma^2$.

In this section, we suppose that
the prior density for $\tau$, conditional on $\sigma$, is $\xi \delta(\tau) + (1-\xi)$,
where $0 \le \xi \le 1$.
In other words, we assume that the prior density of $(\theta, \tau, \sigma^2)$ is
\begin{equation}
\label{first_prior_density_overall}
\pi(\theta, \tau, \sigma^2) = \big(\xi \delta(\tau) + (1-\xi) \big) \sigma^{-2}.
\end{equation}
An attractive feature of this prior density is that, as shown in Appendix G, the $1-\alpha$ HPD and
equi-tailed credible intervals for $\theta$ are identical to the usual frequentist
$1-\alpha$ confidence interval in the two
extreme cases that (a) $\tau$ is known to be 0 (i.e. $\xi = 1$) and (b) there is no prior information about $\tau$
(i.e. $\xi = 0$).

Let $f_q(\, \cdot \,; \mu, \sigma^2)$ denote the density
function of $\mu + \sigma T$, where $\sigma>0$ and $T \sim t_q$. Note that
\begin{equation*}
f_q(x; \mu, \sigma^2) = \frac{1}{\sigma} \, f \left( \frac{x-\mu}{\sigma} \, \bigg | \, t_q\right ),
\end{equation*}
where $f(\, \cdot \, | \, t_q)$ denotes the $t_q$ density function.
Also let $\rho = \text{Corr}(\hat \theta,\hat \tau) = \boldsymbol{a}^T (\boldsymbol{X}^T \boldsymbol{X})^{-1} \boldsymbol{c}$.
Since $\rho$ is determined by $\boldsymbol{a}$, $\boldsymbol{c}$ and $\boldsymbol{X}$, which
are known, we assume that $\rho$ is given.
In Appendix H we consider the prior density \eqref{first_prior_density_overall}.
As shown in this appendix, the marginal posterior density of $\theta$ is
\begin{equation}
\label{first_marg_post_density_theta}
\lambda(\hat{\sigma}, \hat{\tau}/\hat{\sigma})\, f_{m+1}(\theta; \mu_1, {\sigma_1}^2(2))
 + \left(1-\lambda(\hat{\sigma}, \hat{\tau}/\hat{\sigma}) \right) \, f_{m}(\theta; \hat{\theta}, \hat{\sigma}^2),
\end{equation}
where $\mu_1 = \hat{\theta}-\rho \hat{\tau}$, $\sigma_1^2(2) = (m\hat{\sigma}^2 + \hat{\tau}^2)(1-\rho^2)/(m+1)$,
\begin{equation*}
\lambda(\hat{\sigma}, \hat{\tau}/\hat{\sigma}) =  \frac{1}{1 + k \, \hat{\sigma}\,{\left(m + \big(\hat{\tau}/\hat{\sigma}\big)^2 \right)}^{(m+1)/2}}
\end{equation*}
and
\begin{equation*}
k = \frac{(1-\xi)\,\sqrt{\pi}\,\,\Gamma(m/2)}{\xi\,m^{m/2} \, \Gamma((m+1)/2)}.
\end{equation*}
It is easy to find values of $\xi$, $\hat{\sigma}$, $\hat{\tau}/\hat{\sigma}$  and $1-\alpha$ such that the
marginal posterior density of $\theta$ is bimodal and leads to $1-\alpha$ HPD credible
sets that consist of the union of two disjoint intervals.
An illustration is provided by Figure 1.
We therefore focus on Bayesian
$1-\alpha$ equi-tailed and
shortest credible intervals (discussed, for example, by Ferentinos and Karakostas, 2006).
An attractive property of the Bayesian
$1-\alpha$ shortest credible interval is that if the marginal posterior density of $\theta$ is
unimodal then this credible interval is the same as the Bayesian
$1-\alpha$ HPD credible set.
All of the computations presented in this paper were performed with programs
written in MATLAB, using the Optimization and Statistics toolboxes.


\begin{figure}[t]
\centering
\includegraphics[scale=0.8]{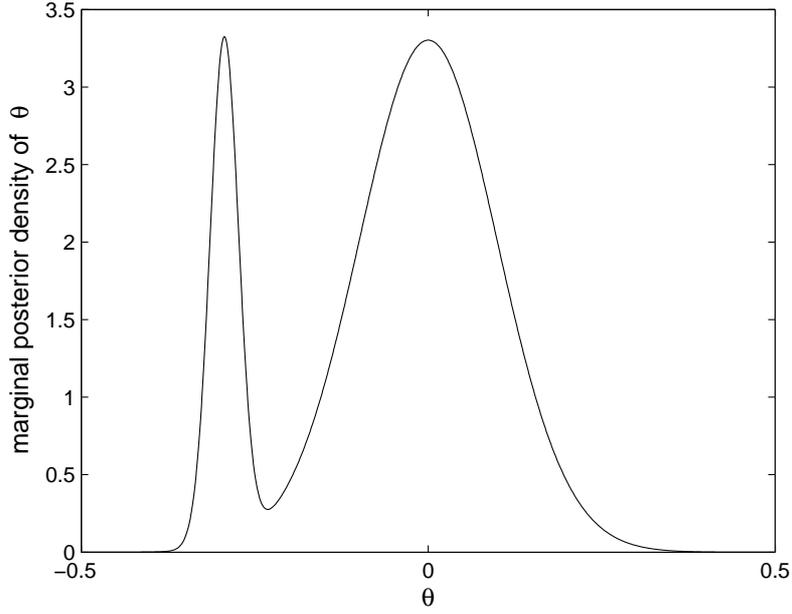}
\caption{\small{Plot of the marginal posterior density of $\theta$, when
the prior density $\pi(\theta, \tau, \sigma^2)$ is $\big(\xi \delta(\tau) + (1-\xi) \big) \sigma^{-2}$, $\xi=0.8$, $\rho=0.98$,
$m=100$, $\hat{\theta}=0$, $\hat{\sigma}=0.1$ and $\hat{\tau}/\hat{\sigma}=3$.}}
\label{Fig1}
\end{figure}


Define $\ell_B(\hat{\theta}, \hat{\tau}, \hat{\sigma}; \eta)$ to be the solution for $\ell$ of
$P(\theta < \ell \, | \, \hat{\theta}, \hat{\tau}, \hat{\sigma}) = \eta$, where
$P(\, \cdot \, | \, \hat{\theta}, \hat{\tau}, \hat{\sigma})$ denotes the posterior probability.
Also define $u_B(\hat{\theta}, \hat{\tau}, \hat{\sigma}; \delta)$ to be the solution for $u$ of
$P(\theta > u \, | \, \hat{\theta}, \hat{\tau}, \hat{\sigma}) = \delta$.
The Bayesian
$1-\alpha$ equi-tailed credible interval is
$\big[\ell_B(\hat{\theta}, \hat{\tau}, \hat{\sigma}; \alpha/2), \, u_B(\hat{\theta}, \hat{\tau}, \hat{\sigma}; \alpha/2) \big]$.
The Bayesian $1-\alpha$ shortest credible interval is
$\big[\ell_B(\hat{\theta}, \hat{\tau}, \hat{\sigma}; \eta^*), \, u_B(\hat{\theta}, \hat{\tau}, \hat{\sigma}; \alpha - \eta^*) \big]$,
where $\eta^*$ minimizes the length of the $1-\alpha$ credible interval
$\big[\ell_B(\hat{\theta}, \hat{\tau}, \hat{\sigma}; \eta), \, u_B(\hat{\theta}, \hat{\tau}, \hat{\sigma}; \alpha - \eta) \big]$
with respect to $\eta \in (0, \alpha)$. Observe that
\begin{align}
\label{sc_half_length_transf_endpoints}
\text{scaled half-length}
&= \frac{1}{2} \left (\frac{u_B(\hat{\theta}, \hat{\tau}, \hat{\sigma}; \delta) - \hat{\theta}}{\hat{\sigma}}
- \frac{\ell_B(\hat{\theta}, \hat{\tau}, \hat{\sigma}; \eta) - \hat{\theta}}{\hat{\sigma}} \right ) \\
\label{sc_offset_transf_endpoints}
\text{scaled offset}
&= \frac{1}{2} \left (\frac{u_B(\hat{\theta}, \hat{\tau}, \hat{\sigma}; \delta) - \hat{\theta}}{\hat{\sigma}}
+ \frac{\ell_B(\hat{\theta}, \hat{\tau}, \hat{\sigma}; \eta) - \hat{\theta}}{\hat{\sigma}} \right ).
\end{align}

Suppose that $\xi \in (0,1)$. As proved in Appendix I, both
$\big(\ell_B(\hat{\theta}, \hat{\tau}, \hat{\sigma}; \eta) - \hat{\theta} \big) / \hat{\sigma}$ and
$\big(u_B(\hat{\theta}, \hat{\tau}, \hat{\sigma}; \delta) - \hat{\theta} \big) / \hat{\sigma}$ are
functions of $(\hat{\sigma}, \hat{\tau} / \hat{\sigma})$. It follows from
this that the scaled half-length and the scaled offset cannot both be functions of $\hat{\tau} / \hat{\sigma}$.
This can be proved by contradiction as follows. Suppose that the scaled half-length and the scaled offset
are both functions of $\hat{\tau} / \hat{\sigma}$.
The sum and difference of the scaled half-length and the scaled offset must also be functions of
$\hat{\tau} / \hat{\sigma}$.
This establishes a contradiction.

Also observe that, irrespective of how large $|\hat{\tau}|/\hat{\sigma}$ is, we can find $\hat{\sigma}$
sufficiently small that the Bayesian $1-\alpha$
equi-tailed and
shortest credible intervals do {\sl not} approximate the interval
$\big[\hat{\theta} - t(m) \hat{\sigma}, \, \hat{\theta} + t(m) \hat{\sigma} \big]$. By contrast,
the interval $J(b,s)$ reverts to the interval
$\big[\hat{\theta} - t(m) \hat{\sigma}, \, \hat{\theta} + t(m) \hat{\sigma} \big]$
when $|\hat{\tau}|/\hat{\sigma} \ge d$. In summary: the Bayesian $1-\alpha$
equi-tailed and shortest credible intervals  depend on the data in
a very different way from this frequentist confidence interval $J(b,s)$.

\medskip

For further numerical illustration, we consider the following example.

\smallskip

\noindent \underbar {\textbf{$\boldsymbol{2 \times 2}$ factorial experiment example}}

\smallskip

\noindent Consider a $2 \times 2$ factorial experiment with
2 replicates and parameter of interest $\theta$
the {\it simple} effect (expected response when factor A is high and factor B is low) $-$
(expected response when factor A is low and factor B is low).
Suppose that we have uncertain prior information that the two-factor
interaction is zero.

Let $x_1$ take the values $-1$ and 1 when the factor A takes the values
low and high respectively. Also let $x_2$ take the values $-1$ and 1 when the factor B takes the values
low and high respectively. In other words, $x_1$ and $x_2$ are the coded values of the
factors A and B, respectively. The model for this experiment is
\begin{equation*}
Y = \beta_0 + \beta_1 x_1 + \beta_2 x_2 + \beta_{12} x_1 x_2 + \varepsilon
\end{equation*}
where $Y$ is the response, $\beta_0$, $\beta_1$, $\beta_2$ and $\beta_{12}$ are unknown parameters and the
$\varepsilon$ for different response measurements are independent and identically $N(0, \sigma^2)$ distributed.
In this case, $n = 8$ and $p = 4$, so that $m = n-p = 4$.
Thus $\theta = 2(\beta_1 - \beta_{12})$. Let $\hat \beta_1$ and $\hat \beta_{12}$ denote the least squares
estimators of $\beta_1$ and $\beta_{12}$ respectively. The least squares estimator of $\theta$ is
$\hat \theta = 2(\hat \beta_1 - \hat \beta_{12})$. Our uncertain prior information is that $\beta_{12}=0$. Note that
$\rho = \text{Corr}(\hat \theta,\hat \tau) = -1/\sqrt{2}$.

\medskip

Figures 2 and 3 illustrate the dependence of the scaled offset and scaled half-length on both $\hat{\sigma}$
and $\hat{\tau} / \hat{\sigma}$ for the Bayesian 0.95 equi-tailed and shortest credible intervals for $\theta$,
in the context of this $2 \times 2$ factorial experiment example,
when the prior density $\pi(\theta, \tau, \sigma^2)$ is $\big(\xi \delta(\tau) + (1-\xi) \big) \sigma^{-2}$
and $\xi = 1/1.2$.

\newpage

\FloatBarrier

\begin{figure}[t]
\centering
\includegraphics[scale=0.7]{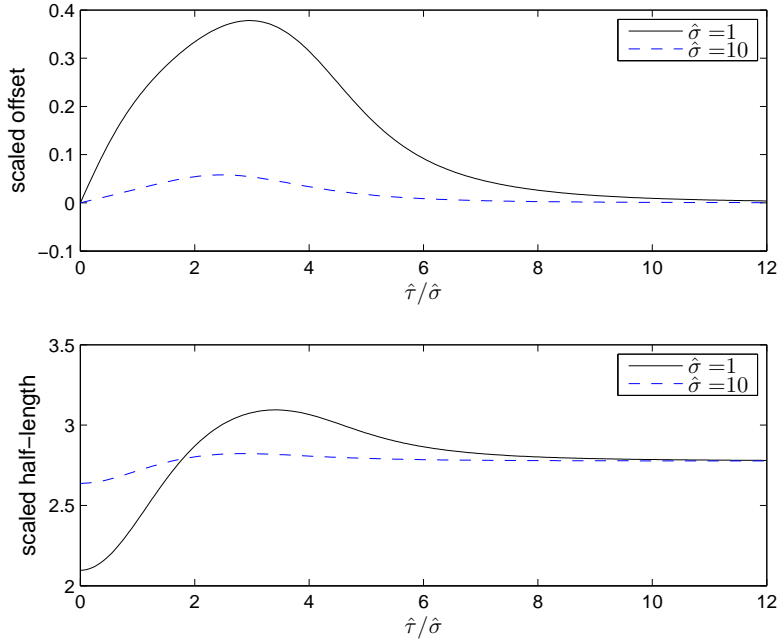}
\caption{\small{Graphs of the scaled offset and scaled half-length, as functions of $\hat{\tau} / \hat{\sigma}$, for $\hat{\sigma} = 1$ (solid line) and $\hat{\sigma} = 10$ (dashed line). These are for the Bayesian 0.95 equi-tailed credible interval for $\theta$,
in the context of the $2 \times 2$ factorial experiment example,
when the prior density $\pi(\theta, \tau, \sigma^2)$ is $\big(\xi \delta(\tau) + (1-\xi) \big) \sigma^{-2}$
and $\xi = 1/1.2$.
}}
\label{Fig2}
\end{figure}

\begin{figure}[h]
\centering
\includegraphics[scale=0.7]{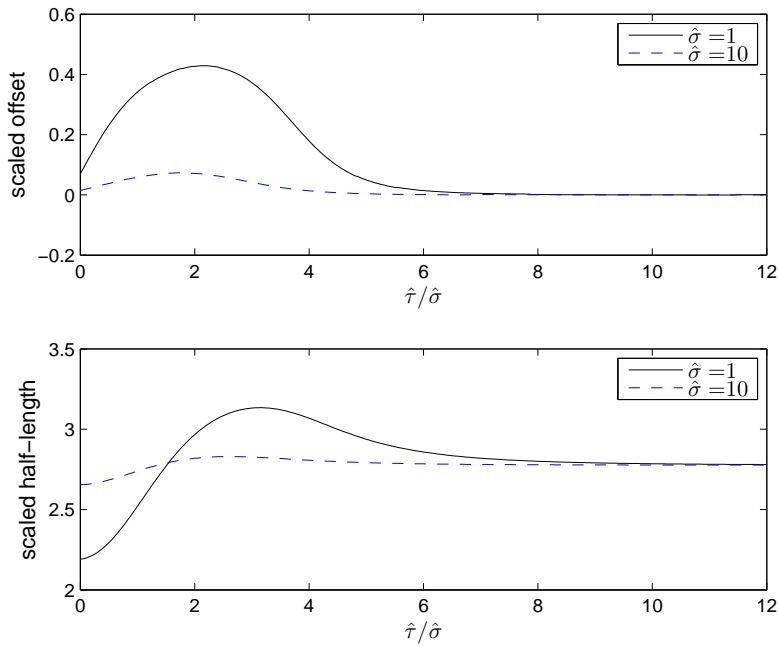}
\caption{\small{Graphs of the scaled offset and scaled half-length, as functions of $\hat{\tau} / \hat{\sigma}$, for $\hat{\sigma} = 1$ (solid line) and $\hat{\sigma} = 10$ (dashed line). These are for the Bayesian 0.95 shortest credible interval for $\theta$,
in the context of the $2 \times 2$ factorial experiment example,
when the prior density $\pi(\theta, \tau, \sigma^2)$ is $\big(\xi \delta(\tau) + (1-\xi) \big) \sigma^{-2}$
and $\xi = 1/1.2$.}}
\label{Fig3}
\end{figure}

\FloatBarrier

\bigskip

\noindent {\large{\bf 4. Bayesian interval estimators for
the prior density \newline $\boldsymbol{\pi(\theta, \tau, \sigma^2)  = \xi \delta(\tau)\sigma^{-1} + (1-\xi)\sigma^{-2}}$}}

\medskip

Let $\gamma = \tau / \sigma$. Since we assume that $\sigma >0$, the uncertain prior information that $\tau=0$ can also
be expressed as the uncertain prior information that $\gamma = 0$.
Suppose that, conditional on $\sigma$, $\gamma$ has the
improper prior density $\xi \delta(\gamma) + (1-\xi)$, where $\xi \in [0,1]$.
Transformations of improper prior densities are problematic. Nonetheless,
the plausibility argument presented in Appendix J suggests that this corresponds to $\tau$ having
prior density $\xi \delta(\tau) + (1-\xi)\sigma^{-1}$, conditional on $\sigma$.
We assume throughout the paper that the prior distributions
of $\theta$ and $(\tau,\sigma^2)$ are independent and that $\theta$ has a uniform prior distribution over the real line.
Assuming that $\sigma^2$ has the standard noninformative prior density $1/\sigma^2$, we obtain the prior density
$\pi(\theta, \tau, \sigma^2)= \xi \delta(\tau)\sigma^{-2} + (1-\xi)\sigma^{-3}$.
Interestingly, for this prior density, both the scaled half-length and the scaled offset are functions of
$\hat{\tau}/\hat{\sigma}$.
In fact, it follows from
\eqref{sc_half_length_transf_endpoints}
and
\eqref{sc_offset_transf_endpoints}
and the results derived
in Appendix L, that this is true for all prior densities of the form
$\pi(\theta, \tau, \sigma^2)  = \xi \delta(\tau){\sigma}^{-g} + (1-\xi) {\sigma}^{-g-1}$,
where $m+g > 0$.

We focus on the particular case that $g=1$, so that the prior density is
 $\pi(\theta, \tau, \sigma^2)= \xi \delta(\tau)\sigma^{-1} + (1-\xi)\sigma^{-2}$.
 We have chosen to focus on this prior density because, for $\xi = 0$ i.e. for
no prior information about $\tau$,
the $1-\alpha$ HPD and equi-tailed credible intervals for $\theta$ are identical
to the usual frequentist $1-\alpha$ confidence interval $I$.
As shown in Appendix K, the posterior marginal density of $\theta$ is equal to
\begin{equation*}
\tilde{\lambda}(\hat{\tau}/\hat{\sigma},1) \, f_{m}(\theta; \mu_1, \sigma_1^2(1))
+(1-\tilde{\lambda}(\hat{\tau}/\hat{\sigma},1)) \, f_{m}(\theta; \hat{\theta}, \hat{\sigma}^2),
\end{equation*}
where $\mu_1 = \hat{\theta}-\rho\hat{\tau}$, $\sigma_1^2(1)=(m\hat{\sigma}^2+\hat{\tau}^2)(1-\rho^2)/m$,
\begin{equation*}
\tilde{\lambda}(\hat{\tau}/\hat{\sigma},1)= \frac{1}{1+\tilde{k}(1)\big(m+(\hat{\tau}/\hat{\sigma})^2\big)^{(m+1)/2}}
\end{equation*}
and
\begin{equation*}
\tilde{k}(1) =\frac {\sqrt{2\pi}\,(1-\xi)}{\xi m^{m/2}}.
\end{equation*}
For the case that $\tau$ is known to be 0 (i.e. $\xi=1$),
the posterior marginal density for $\theta$ is $f_{m}(\theta; \mu_1, \sigma_1^2(1))$, instead of the density
$f_{m+1} \big(\theta; \mu_1, \sigma_1^2(2) \big)$ which results from the prior
density $\pi(\theta, \tau, \sigma^2)  = \big(\xi \delta(\tau) + (1-\xi) \big) {\sigma}^{-2}$.
Therefore, for this case, the $1-\alpha$ HPD and equi-tailed credible intervals for $\theta$ are {\sl not} the same as
the usual frequentist $1-\alpha$ confidence interval for $\theta$, assuming that $\tau=0$.
On the other hand,
the $1-\alpha$ credible intervals for $\theta$ based on
$f_{m}(\theta; \mu_1, \sigma_1^2(1))$ and $f_{m+1} \big(\theta; \mu_1, \sigma_1^2(2) \big)$ are approximately equal for
large $m$ in the sense that they are both centered on $\mu_1$ and the ratio of their lengths approaches 1,
 as $m \rightarrow \infty$.

 Also observe that $\tilde{\lambda}(\hat{\tau}/\hat{\sigma})\rightarrow 0$ as
$|\hat{\tau}|/\hat{\sigma} \rightarrow \infty$. Therefore,
 the Bayesian $1-\alpha$ equi-tailed and
shortest credible intervals approach the interval
$\big[\hat{\theta} - t(m) \hat{\sigma}, \, \hat{\theta} + t(m) \hat{\sigma} \big]$
when $|\hat{\tau}|/\hat{\sigma}$ is large.
Similarly, the interval $J(b,s)$ reverts to this interval
when $|\hat{\tau}|/\hat{\sigma} \ge d$.

\medskip

Figures 4 and 5 illustrate the fact that the scaled offset and scaled half-length are functions of
$\hat{\tau} / \hat{\sigma}$ for the Bayesian 0.95 equi-tailed and shortest credible intervals for $\theta$,
in the context of the $2 \times 2$ factorial experiment example,
when the prior density $\pi(\theta, \tau, \sigma^2)$ is $\xi \delta(\tau)\sigma^{-1} + (1-\xi)\sigma^{-2}$
and $\xi = 1/1.2$.

 \newpage

\FloatBarrier

\begin{figure}[t]
\centering
\includegraphics[scale=0.75]{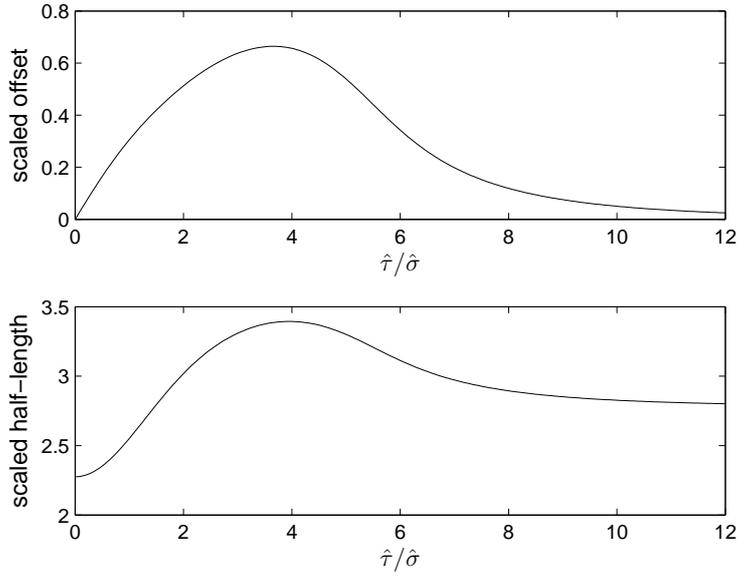}
\caption{\small{Graphs of the scaled offset and the scaled half-length, as functions of $\hat{\tau} / \hat{\sigma}$,
 for the Bayesian 0.95  equi-tailed credible interval
for $\theta$,
in the context of the $2 \times 2$ factorial experiment example, when
the prior density $\pi(\theta, \tau, \sigma^2)$
is $\xi \delta(\tau)\sigma^{-1} + (1-\xi)\sigma^{-2}$ and $\xi = 1/1.2$.}}
\label{Fig4}
\end{figure}

\begin{figure}[h]
\centering
\includegraphics[scale=0.75]{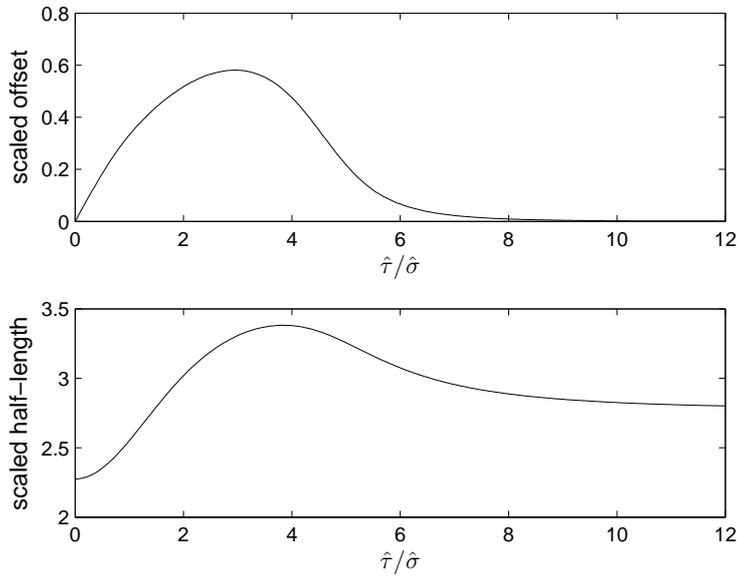}
\caption{\small{
Graphs of the scaled offset and the scaled half-length, as functions of $\hat{\tau} / \hat{\sigma}$,
 for the Bayesian 0.95  shortest credible interval
for $\theta$,
in the context of the $2 \times 2$ factorial experiment example, when
the prior density $\pi(\theta, \tau, \sigma^2)$
is $\xi \delta(\tau)\sigma^{-1} + (1-\xi)\sigma^{-2}$ and $\xi = 1/1.2$.}}
\label{Fig5}
\end{figure}

\FloatBarrier

\bigskip

\noindent {\large{\bf 5. Comparison with the KG
$\boldsymbol{1-\alpha}$ confidence interval}}

\medskip

We begin this section by briefly describing a method for computing the KG $1-\alpha$ confidence interval for
$\theta$ that utilizes the uncertain prior information that $\tau=0$, in the sense described in the introduction.
For given $(b,s)$,
the coverage probability $P\big( \theta \in J(b, s) \big)$ is an even function of $\gamma$.
The scaled expected length of $J(b,s)$ is (expected length of $J(b, s)$)/(expected length of $I$)
and is an even function of $\gamma$ for given $s$, which we denote
by $e(\gamma;s)$.
Define the weight function
$w(\gamma) = \tilde{\xi} \delta(\gamma) + (1-\tilde{\xi})$,
where $0 \le \tilde{\xi} \le 1$.
Kabaila and Giri (2009) describe how to compute smooth functions $b$ and $s$
such that (a) the minimum of $P\big( \theta \in J(b, s) \big)$ over $\gamma$ is
$1-\alpha$ and (b)
\begin{equation}
\label{criterion}
\int_{-\infty}^{\infty} (e(\gamma;s) - 1) \, w(\gamma) \, d \gamma
\end{equation}
is minimized, where $\tilde{\xi}$ is a specified tuning parameter.
This tuning parameter and the functions $b$ and $s$ are chosen by the statistician
{\sl prior} to looking at the observed
response vector $\boldsymbol{y}$.

Consider the $2 \times 2$ factorial experiment example.
For $d=6$, functions $b$ and $s$ chosen to be natural cubic splines in the interval $[0,d]$,
with evenly-spaced knots at $0, 2, 4, \dots , 12$, $1-\alpha = 0.95$ and $\tilde{\xi} = 1/1.2$, the minimization
of \eqref{criterion}, subject to the coverage constraint, leads to a confidence interval with
the following properties. To within computational accuracy, this confidence interval has coverage probability
0.95 for all $\gamma$ (i.e. throughout the parameter space). Figure 6 is a plot of the squared scaled expected length
of this confidence interval, as a function of $\gamma$.
It is clear from this figure that this confidence interval utilizes the uncertain prior information that $\tau=0$, in the sense described in the introduction.
 When
the prior information is correct (i.e. $\gamma=0$), we gain since
$e^2(0;s) = 0.8524$. The maximum value of $e^2(\gamma;s)$ is
1.0852. This confidence interval coincides with the standard
$1-\alpha$ confidence interval for $\theta$ when the data strongly
contradicts the prior information, so that $e^2(\gamma;s)$
approaches 1 as $\gamma \rightarrow \infty$.

Figure 7 shows graphs of the scaled offset and scaled half-length for this confidence interval,
as functions of $\hat{\tau} / \hat{\sigma}$. The weight function
$w(\gamma)$ has a similar form to the prior density of $\gamma = \tau / \sigma$, conditional on $\sigma$,
considered in Section 4, when $\tilde{\xi} = \xi$. Although this weight and conditional prior density
have very different interpretations, it is of interest to compare the scaled offset and scaled half-length
shown in Figures 4 and 4 with the scaled offset and scaled half-length shown in Figure 7.
The differences are marked, particularly with respect to the scaled offset.

\FloatBarrier

\begin{figure}[h]
\centering
\includegraphics[scale=0.75]{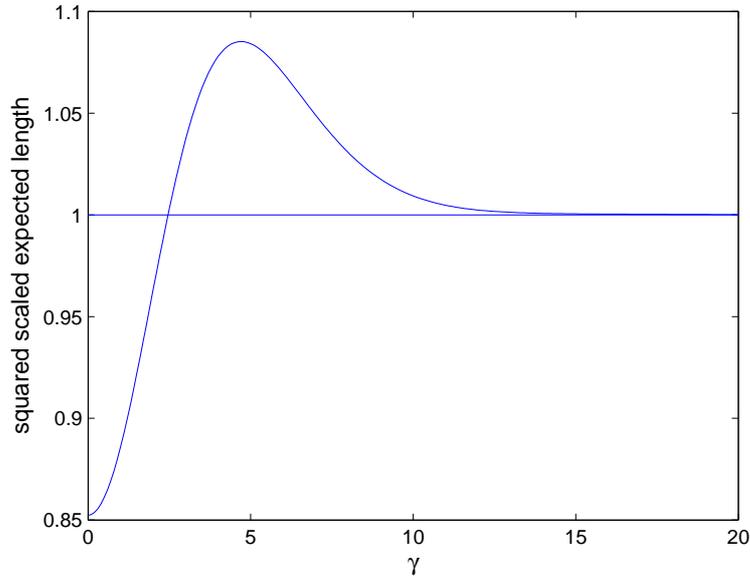}
\caption{\small{
Graph of the squared scaled expected length, as a function of $\gamma$,
 for the KG 0.95  confidence interval
for $\theta$,
in the context of the $2 \times 2$ factorial experiment example, when
$d=12$, $b$ and $s$ are natural cubic splines in the interval $[0,d]$,
with evenly-spaced knots at $0, 2, 4, \dots , 12$ and $\tilde{\xi} = 1/1.2$.}}
\label{Fig6}
\end{figure}

\begin{figure}[h]
\centering
\includegraphics[scale=0.75]{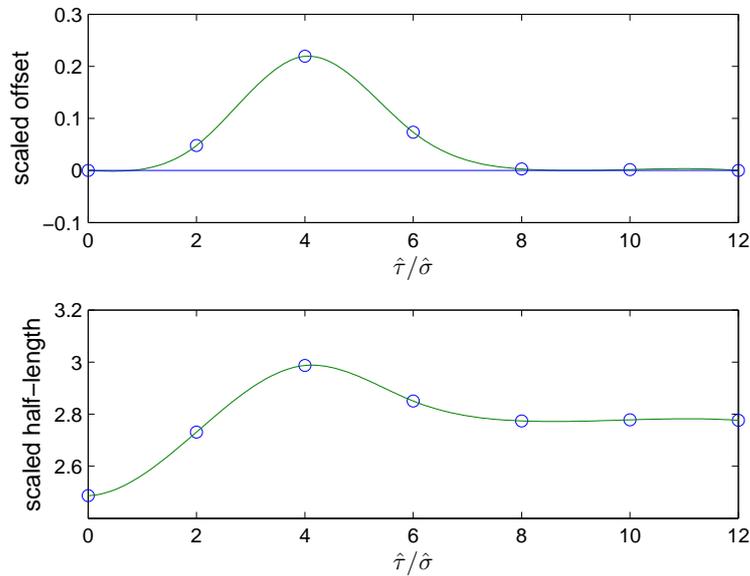}
\caption{\small{
Graphs of the scaled offset and the scaled half-length, as functions of $\hat{\tau} / \hat{\sigma}$,
 for the KG 0.95  confidence interval
for $\theta$,
in the context of the $2 \times 2$ factorial experiment example, when
$d=12$, $b$ and $s$ are natural cubic splines in the interval $[0,d]$,
with evenly-spaced knots at $0, 2, 4, \dots , 12$ and $\tilde{\xi} = 1/1.2$.
The knots of the cubic splines are denoted by small circles.}}
\label{Fig7}
\end{figure}

\FloatBarrier

\bigskip

\noindent {\large{\bf 6. Discussion}}

\medskip

We have not sought to advocate the use of Bayesian $1-\alpha$ credible intervals in place of
frequentist $1-\alpha$ confidence intervals or vice versa.
Bayesian $1-\alpha$ credible intervals could be examined from the point of view of their frequentist
coverage properties. Similarly, frequentist $1-\alpha$ confidence intervals could be examined
from the point of their Bayesian posterior coverage properties. The first comparison is likely
to favor the frequentist $1-\alpha$ confidence intervals, since the Bayesian $1-\alpha$ credible intervals
are not constructed to have good frequentist coverage properties. Similarly, the second comparison
is likely to favor the Bayesian $1-\alpha$ credible intervals. By contrasting the dependencies on the
data of the frequentist $1-\alpha$ confidence intervals and the Bayesian $1-\alpha$ credible intervals,
we have avoided a comparison that is partisan to either frequentist or Bayesian points of view.

Bayesian and frequentist statistical analyses differ in important ways. However, it is pleasing when they lead
to the same result. In the present paper we have found yet another instance of a difference between
Bayesian and frequentist statistical analyses.

\bigskip

\noindent {\large{\bf {Appendix A: Initial scaling of the parameters}}}

\medskip

We assume that Var$(\hat{\theta}) = \sigma^2$ and
Var$(\hat{\tau}) = \sigma^2$. In this Appendix, we show that this can be achieved by appropriate scaling
and that there is no loss of generality, as far as the purposes of the paper are concerned.

Suppose that the parameter of interest is $\theta^* = (\boldsymbol{a}^*)^T \boldsymbol{\beta}$
where $\boldsymbol{a}^*$ is a specified
$p$-vector ($\boldsymbol{a}^* \ne \boldsymbol{0}$).
Suppose that the inference of interest is an interval estimator for $\theta^*$.
Define the parameter $\tau^* = (\boldsymbol{c}^*)^T \boldsymbol{\beta} - t^*$ where the vector $\boldsymbol{c}^*$ and the number $t^*$
are specified and $\boldsymbol{a}^*$ and $\boldsymbol{c}^*$ are linearly independent.
Also suppose that previous experience with similar data sets and/or
expert opinion and scientific background suggest that $\tau^* = 0$.
In other words, suppose that we have uncertain prior information that $\tau^* = 0$.
Let $v_{11} = (\boldsymbol{a}^*)^T (\boldsymbol{X}^T \boldsymbol{X})^{-1} (\boldsymbol{a}^*)$ and
$v_{22} = (\boldsymbol{c}^*)^T (\boldsymbol{X}^T \boldsymbol{X})^{-1} (\boldsymbol{c}^*)$.
It is convenient to transform $\theta^*$ to $\theta = \theta^*/\sqrt{v_{11}} = \boldsymbol{a}^T \boldsymbol{\beta}$,
where $\boldsymbol{a} = \boldsymbol{a}^*/\sqrt{v_{11}}$.
It is also convenient to transform $\tau^*$ to $\tau = \tau^*/\sqrt{v_{22}} = \boldsymbol{c}^T \boldsymbol{\beta} - t$,
where $\boldsymbol{c} = \boldsymbol{c}^*/\sqrt{v_{22}}$ and $t = t^*/\sqrt{v_{22}}$.
Since $\hat{\theta} = \boldsymbol{a}^T \hat{\boldsymbol{\beta}}$ and
$\hat{\tau} = \boldsymbol{c}^T \hat{\boldsymbol{\beta}} - t$, where
$\hat{\boldsymbol{\beta}} \sim N(\boldsymbol{\beta}, \sigma^2 (\boldsymbol{X}^T \boldsymbol{X})^{-1})$,
Var$(\hat{\theta}) = \sigma^2$ and
Var$(\hat{\tau}) = \sigma^2$.
Interval estimators for $\theta$ and their properties transform in the obvious way to interval estimators
for $\theta^*$ and their corresponding properties. The uncertain prior information that $\tau^* = 0$ implies the uncertain prior
information that $\tau = 0$.

\bigskip

\noindent {\large{\bf {Appendix B: Transformation of the regression model}}}

\medskip

Define $\boldsymbol{W} = (\boldsymbol{X}^T \boldsymbol{X})^{-1}.$ Now define the $p \times p$ matrix
$\boldsymbol{B}$ as follows. The first and second rows of $\boldsymbol{B}$ are $\boldsymbol{a}^T {\boldsymbol{W}}^{1/2}$
and $\boldsymbol{c}^T {\boldsymbol{W}}^{1/2}$, respectively. The last $p-2$ rows consist of unit-length
orthogonal $p$ -vectors, that are orthogonal to both $\boldsymbol{a}$ and $\boldsymbol{c}$.

We re-express the regression sampling model as
$\boldsymbol{Y} = \tilde{\boldsymbol{X}} \boldsymbol{\lambda} + \boldsymbol{\varepsilon}$, where\\
$\tilde{\boldsymbol{X}} = \boldsymbol{X} {\boldsymbol{W}}^{1/2} {\boldsymbol{B}}^{-1}$ and
$\boldsymbol{\lambda} = \boldsymbol{B} {\boldsymbol{W}}^{-1/2} \boldsymbol{\beta}$. Let $\hat{\boldsymbol{\lambda}}$
denote the least squares estimator of $\boldsymbol{\lambda}$, based on this model. We reduce the data to the sufficient
statistic $(\hat{\boldsymbol{\lambda}}, \hat{\sigma}^2)$ for $(\boldsymbol{\lambda}, \sigma^2)$. Let
\begin{equation*}
\boldsymbol{V} = \left[\begin{matrix} 1 \quad \rho\\ \rho \quad 1 \end{matrix}
\right].
\end{equation*}
Note that
$\hat{\boldsymbol{\lambda}} \sim N \big(\boldsymbol{\lambda}, \sigma^2 \boldsymbol{D}\big)$, where
\begin{equation*}
\boldsymbol{D} =
\left[
\begin{array}{c|c}
 \boldsymbol{V} & \boldsymbol{0}  \\ \hline \\[-9pt]
  \boldsymbol{0} & \boldsymbol{I}_{p-2} \\[-1pt]
 \end{array}
\right].
\end{equation*}
Observe that $\theta = \lambda_1$ and
$\tau = \lambda_2 - t$. Define the parameter vector
$\boldsymbol{\chi} = (\lambda_3,\dots, \lambda_p)$. Now define
$\hat{\boldsymbol{\chi}} = (\hat{\lambda}_3,\dots, \hat{\lambda}_p)$.
Since $\boldsymbol{D}$ is block diagonal,
$(\hat{\theta}, \hat{\tau})$ and $\hat{\boldsymbol{\chi}}$ are independent random vectors,
$(\hat{\theta}, \hat{\tau})$ has the distribution \eqref{model_thetahat_tauhat} and
$\hat{\boldsymbol{\chi}} \sim N \big(\boldsymbol{\chi}, \sigma^2 \boldsymbol{I}_{p-2}\big)$.

\bigskip

\noindent {\large{\bf {Appendix C: Marginal posterior distribution of $\boldsymbol{(\theta, \tau, \sigma^2)}$}}}

\medskip

Suppose that the prior distributions of
$(\theta, \tau, \sigma^2)$ and $\boldsymbol{\chi}$  are independent. Also suppose that
the components of $\boldsymbol{\chi}$  have independent uniform prior distributions. In this appendix we prove that
the marginal posterior distribution of $(\theta, \tau, \sigma^2)$ is the same as the posterior
distribution of $(\theta, \tau, \sigma^2)$ based on the reduced data
$\big(\hat{\theta}, \hat{\tau}, \hat \sigma^2 \big)$ and the sampling model that
$\big(\hat{\theta}, \hat{\tau}\big)$ and $\hat \sigma^2$ are independent random vectors, $\big(\hat{\theta}, \hat{\tau}\big)$
has the distribution \eqref{model_thetahat_tauhat} and $(n-p) \hat{\sigma}^2/ \sigma^2 \sim \chi^2_{n-p}$.

It follows from Appendix B that, under the sampling model,
\begin{equation*}
f(\hat{\theta}, \hat{\tau}, \hat{\boldsymbol{\chi}}, \hat{\sigma}^2|\theta, \tau, \boldsymbol{\psi}, \sigma^2) =
f(\hat{\theta}, \hat{\tau}|\theta, \tau, \sigma^2) f(\hat{\boldsymbol{\chi}}|\boldsymbol{\chi}, \sigma^2) f(\hat{\sigma}^2|\sigma^2)
\end{equation*}
where $f(\hat{\boldsymbol{\chi}}|\boldsymbol{\chi}, \sigma^2) = (2\pi \sigma^2)^{-(p-2)/2} \, \exp
\big(-{(\hat{\boldsymbol{\chi}}-\boldsymbol{\chi})}^T(\hat{\boldsymbol{\chi}}-\boldsymbol{\chi})/2\sigma^2\big)$.

Suppose that the prior distributions of $(\theta, \tau, \sigma^2)$ and $\boldsymbol{\chi}$ are independent.
Let $\pi(\theta, \tau, \sigma^2)$ denote the prior density of $(\theta, \tau, \sigma^2)$. Suppose that the components
of $\boldsymbol{\chi}$ are independent and uniformly distributed. Thus the prior density of $(\theta, \tau, \boldsymbol{\chi}, \sigma^2)$
is $\pi(\theta, \tau, \sigma^2)$. Hence the posterior density $\pi(\theta, \tau, \boldsymbol{\chi}, \sigma^2|\hat{\theta}, \hat{\tau}, \hat{\boldsymbol{\chi}}, \hat{\sigma}^2)$
is proportional to $f(\hat{\theta}, \hat{\tau}|\theta, \tau, \sigma^2) f(\hat{\boldsymbol{\chi}}|\boldsymbol{\chi}, \sigma^2) f(\hat{\sigma}^2|\sigma^2)
\pi(\theta, \tau, \sigma^2)$. Thus the marginal posterior density of $(\theta, \tau, \sigma^2)$ is proportional to
\begin{equation*}
 f(\hat{\theta}, \hat{\tau}|\theta, \tau, \sigma^2) \boldsymbol{\int} f(\hat{\boldsymbol{\chi}}|\boldsymbol{\chi}, \sigma^2) d\boldsymbol{\chi} \, f(\hat{\sigma}^2|\sigma^2) \,
\pi(\theta, \tau, \sigma^2) = f(\hat{\theta}, \hat{\tau}|\theta, \tau, \sigma^2) f(\hat{\sigma}^2|\sigma^2) \pi(\theta, \tau, \sigma^2).
\end{equation*}

\bigskip

\noindent {\large{\bf {Appendix D: Two useful integrals}}}

\medskip
We make extensive use of the following two integrals:
\begin{equation}
\label{first_integral}
\int_{0}^{\infty} x^{-(p+1)} \, \exp \left(-ax^{-1} \right) \, dx = a^{-p} \, \Gamma(p),
\end{equation}
for $a > 0$ and $p > 0$, which is (A2.1.2) on p.144 of Box and Tiao(1973) and
can be proved by a change of variable of integration in the definition of the gamma function, and
\begin{equation}
\label{second_integral}
\int_{-\infty}^{\infty} \exp \left(-\frac{x^2-2\rho xy+y^2}{2\sigma^2(1-\rho^2)}\right) \, dy = \sqrt{2\pi (1-\rho^2)} \, \sigma \,
\exp \left(-\frac{x^2}{2\sigma^2}\right)
\end{equation}
which can be proved by completion of the square. Note that  \eqref{second_integral} is used in the derivation of a marginal
density for a bivariate normal distribution.

\bigskip

\noindent {\large{\bf {Appendix E: Marginal posterior distribution of $\boldsymbol{\theta}$ for the
prior density $\boldsymbol{\pi(\theta, \tau, \sigma^2)= \delta(\tau)\sigma^{-g}}$}}}

\medskip

In this appendix we suppose that the prior density is $\pi(\theta, \tau, \sigma^2) = \delta(\tau) {\sigma}^{-g}$,
where $m + g > 0$. We derive the marginal
posterior density of $\theta$.
The likelihood function is proportional to
$\tilde{\ell} \big(\theta, \tau, \sigma^2 \, | \, \hat{\theta},\hat{\tau},\hat{\sigma}^2 \big)$,
which is defined to be
\begin{equation}
\label{lik_prop}
\sigma^{-(m+2)} \exp \left(-\frac{m \hat{\sigma}^2}{2\sigma^2} \right)
\exp \left(-\frac{(\theta-\hat{\theta})^2-2\rho(\theta-\hat{\theta})(\tau-\hat{\tau})+(\tau-\hat{\tau})^2}{2\sigma^2(1-\rho^2)}\right).
\end{equation}
Thus the posterior density of $(\theta, \tau, \sigma^2)$
is proportional to
$\tilde{\ell} \big(\theta, \tau, \sigma^2 \, | \, \hat{\theta},\hat{\tau},\hat{\sigma}^2 \big) \delta(\tau) {\sigma}^{-g}$,
which is equal to
\begin{equation}
\label{post_prop_tau_0}
\sigma^{-(m+g+2)} \exp \left(-\frac{m \hat{\sigma}^2}{2\sigma^2} \right)
\exp \left(-\frac{(\theta-\hat{\theta})^2-2\rho(\theta-\hat{\theta})(\tau-\hat{\tau})+(\tau-\hat{\tau})^2}{2\sigma^2(1-\rho^2)}\right) \delta(\tau).
\end{equation}
The marginal posterior density of $\theta$ is proportional to
\begin{equation}
\label{marg_post_density_prop_to}
\int_0^{\infty} \int_{-\infty}^{\infty}
\tilde{\ell} \big(\theta, \tau, \sigma^2 \, | \, \hat{\theta},\hat{\tau},\hat{\sigma}^2 \big) \delta(\tau) {\sigma}^{-g}
\, d \tau \, d \sigma^2
\end{equation}
which is found as follows. Observe that
\begin{align}
\label{Application_Dirac delta}
&\int_{-\infty}^{\infty} \exp \left(-\frac{(\theta-\hat{\theta})^2-2\rho(\theta-\hat{\theta})(\tau-\hat{\tau})+(\tau-\hat{\tau})^2}{2\sigma^2(1-\rho^2)}\right)
\delta(\tau) d\tau \notag\\
&= \exp \left(-\frac{(\theta-\hat{\theta})^2+2\rho(\theta-\hat{\theta})\hat{\tau}+\hat{\tau}^2}{2\sigma^2(1-\rho^2)}\right) \notag\\
&= \exp \left(-\frac{\big(\theta-(\hat{\theta}-\rho \hat{\tau})\big)^2}{2\sigma^2(1-\rho^2)} - \frac{\hat{\tau}^2}{2\sigma^2}\right).
\end{align}
Thus, \eqref{marg_post_density_prop_to} is equal to
\begin{align}
\label{marg_post_prop_tau_0}
 &\int_{0}^{\infty} \sigma^{-(m+g+2)} \exp \left( -\left(\frac{m \hat{\sigma}^2+\hat{\tau}^2}{2} +\frac{\big(\theta-(\hat{\theta}-\rho \hat{\tau})\big)^2}{2(1-\rho^2)} \right)
 {\sigma}^{-2}\right) d\sigma^2  \notag\\
 &=\Gamma\left(\frac{m+g}{2}\right) \left(\frac{m \hat{\sigma}^2+\hat{\tau}^2}{2} +\frac{\big(\theta-(\hat{\theta}-\rho \hat{\tau})\big)^2}{2(1-\rho^2)}\right)^{-(m+g)/2}
\end{align}
by \eqref{first_integral}. Now \eqref{marg_post_prop_tau_0} can be shown to be equal to
$c_1(\hat{\tau}, \hat{\sigma}^2,g) f_{m+g-1}(\theta; \mu_1, {\sigma_1}^2(g))$,
where
\begin{equation}
\label{c1}
c_1(\hat{\tau}, \hat{\sigma}^2, g) = \sqrt{2\pi (1-\rho^2)} {\left(\frac{m \hat{\sigma}^2+\hat{\tau}^2}{2}\right)}^{-(m+g-1)/2} \Gamma\left(\frac{m+g-1}{2}\right)
\end{equation}
and $\mu_1 = \hat{\theta}-\rho \hat{\tau}$ and $\sigma_1^2(g) = (m\hat{\sigma}^2 + \hat{\tau}^2)(1-\rho^2)/(m+g-1)$.
Thus the marginal posterior density of $\theta$ is
$f_{m+g-1}(\theta; \mu_1, {\sigma_1}^2(g))$, where $f_q(\, \cdot \,; \mu, \sigma^2)$ denotes the density
function of $\mu + \sigma T$, where $\sigma>0$ and $T \sim t_q$. Note that
\begin{equation*}
f_q(x; \mu, \sigma^2) = \frac{1}{\sigma} \, f \left( \frac{x-\mu}{\sigma} \, \bigg | \, t_q\right ),
\end{equation*}
where $f(\, \cdot \, | \, t_q)$ denotes the $t_q$ density function. For the particular case that $g=2$, the marginal prior density
of $\theta$ is $f_{m+1}(\theta; \mu_1, {\sigma_1}^2(2))$.

\bigskip

\noindent {\large{\bf {Appendix F: Marginal posterior distribution of $\boldsymbol{\theta}$ for the
prior density $\boldsymbol{\pi(\theta, \tau, \sigma^2)= \sigma^{-h}}$}}}

\medskip

In this appendix we suppose that the prior density is $\pi(\theta, \tau, \sigma^2) = {\sigma}^{-h}$,
where $m+h-1 > 0$. We derive the marginal
posterior density of $\theta$. As noted in Appendix E, the likelihood function is proportional to
$\tilde{\ell} \big(\theta,\tau,\sigma^2 \, \big| \, \hat{\theta},\hat{\tau},\hat{\sigma}^2 \big)$, defined to be \eqref{lik_prop}.
Thus the posterior density of $(\theta,\tau,\sigma^2)$ is proportional to
$\tilde{\ell}\big(\theta,\tau,\sigma^2 \, \big| \, \hat{\theta},\hat{\tau},\hat{\sigma}^2 \big) {\sigma}^{-h}$,
which is equal to
\begin{equation}
\label{post_prop_uniform_tau}
\sigma^{-(m+h+2)} \exp \left(-\frac{m \hat{\sigma}^2}{2\sigma^2} \right)
\exp \left(-\frac{(\theta-\hat{\theta})^2-2\rho(\theta-\hat{\theta})(\tau-\hat{\tau})+(\tau-\hat{\tau})^2}{2\sigma^2(1-\rho^2)}\right).
\end{equation}
The marginal posterior density of $\theta$ is proportional to
\begin{equation}
\label{marg_post_density_prop_to_2}
\int_0^{\infty} \int_{-\infty}^{\infty}
\tilde{\ell} \big(\theta, \tau, \sigma^2 \, | \, \hat{\theta},\hat{\tau},\hat{\sigma}^2 \big) {\sigma}^{-h}
\, d \tau \, d \sigma^2
\end{equation}
which is found
as follows. Note that, by \eqref{second_integral},
\begin{equation*}
\int_{-\infty}^{\infty} \exp \left(-\frac{(\theta-\hat{\theta})^2-2\rho(\theta-\hat{\theta})(\tau-\hat{\tau})+(\tau-\hat{\tau})^2}{2\sigma^2(1-\rho^2)}\right) d\tau
\end{equation*}
is equal to $\sqrt{2\pi (1-\rho^2)} \, \sigma \exp \left(-(\theta-\hat{\theta})^2/2\sigma^2\right)$.
Thus \eqref{marg_post_density_prop_to_2} is equal to
\begin{align}
\label{marg_post_density}
 &\sqrt{2\pi (1-\rho^2)}\int_{0}^{\infty} \sigma^{-(m+h+1)} \exp \left( -\left(\frac{m \hat{\sigma}^2+(\theta-\hat{\theta})^2}{2} \right) {\sigma}^{-2}\right) d\sigma^2
 \notag\\
 &=\sqrt{2\pi (1-\rho^2)} \, \Gamma\left(\frac{m+h-1}{2}\right) \left(\frac{m \hat{\sigma}^2+(\theta-\hat{\theta})^2}{2} \right)^{-(m+h-1)/2}
\end{align}
by \eqref{first_integral}.
Now \eqref{marg_post_density} is equal to $c_2(\hat{\sigma}^2,h) f_{m+h-2}(\theta; \hat{\theta}, {\sigma_2}^2(h))$,
where
\begin{equation}
\label{c_2}
c_2(\hat{\sigma}^2,h) = 2\pi\sqrt{1-\rho^2} {\left(\frac{m \hat{\sigma}^2}{2}\right)}^{-(m+h-2)/2} \Gamma\left(\frac{m+h-2}{2}\right),
\end{equation}
$\sigma_2^2(h) = (m/(m+h-2)) \, \hat{\sigma}^2$. Thus the marginal posterior density of $\theta$ is
$f_{m+h-2}(\theta; \hat{\theta}, {\sigma_2}^2(h))$ For the particular case that $h=2$, this marginal
posterior density is $f_{m}(\theta; \hat{\theta}, \hat{\sigma}^2)$.

\bigskip

\noindent {\large{\bf {Appendix G: An attractive feature of the prior density \newline
$\boldsymbol{\pi(\theta, \tau, \sigma^2)  = \big(\xi \delta(\tau) + (1-\xi) \big) {\sigma}^{-2}}$}}}

\medskip

At the end of Appendix E, we considered the extreme case that $\tau$ is known to be 0  i.e.
$\pi(\theta, \tau, \sigma^2) = \delta(\tau) \sigma^{-2}$. This corresponds to choosing $\xi = 1$
in the prior density $\pi(\theta, \tau, \sigma^2)  = \big(\xi \delta(\tau) + (1-\xi) \big) {\sigma}^{-2}$.
As shown in this appendix,
the marginal posterior density of $\theta$ is
$\pi \big(\theta \, | \, \hat{\theta}, \hat{\tau}, \hat{\sigma}^2 \big) = f_{m+1} \big(\theta; \mu_1, \sigma_1^2(2) \big)$,
where $\mu_1 = \hat{\theta}-\rho \hat{\tau}$ and $\sigma_1^2(2) = (m\hat{\sigma}^2 + \hat{\tau}^2)(1-\rho^2)/(m+1)$.
In this case, the HPD and equi-tailed $1-\alpha$ credible intervals for $\theta$ are identical to
%
%
the usual frequentist $1-\alpha$ confidence interval for $\theta$, assuming that $\tau = 0$.

At the end of Appendix F, we considered the second extreme case that there is no prior information about $\tau$
i.e. $\pi(\theta, \tau, \sigma^2) = \sigma^{-2}$. This corresponds to choosing $\xi = 0$
in the prior density $\pi(\theta, \tau, \sigma^2)  = \big(\xi \delta(\tau) + (1-\xi) \big) {\sigma}^{-2}$.
As shown in this appendix, the marginal posterior
density of $\theta$ is
$\pi(\theta \, | \, \hat{\theta}, \hat{\tau}, \hat{\sigma}^2) = f_{m}(\theta; \hat{\theta}, \hat{\sigma}^2 )$.
In this case also the HPD and equi-tailed $1-\alpha$ credible intervals for $\theta$
are identical and are equal to
$\big[\hat{\theta} - t(m) \hat{\sigma}, \, \hat{\theta} + t(m) \hat{\sigma} \big]$.
This is the same as the usual frequentist $1-\alpha$ confidence interval for $\theta$, assuming that there is no
prior information about $\tau$.

\bigskip

\noindent {\large{\bf {Appendix H: Marginal posterior distribution of $\boldsymbol{\theta}$ for the
prior density $\boldsymbol{\pi(\theta, \tau, \sigma^2)  = \big(\xi \delta(\tau) + (1-\xi) \big) {\sigma}^{-2}}$}}}

\medskip

Suppose that the prior density is $\pi(\theta, \tau, \sigma^2)  = \big(\xi \delta(\tau) + (1-\xi) \big) {\sigma}^{-2}$.
Also suppose that $\xi \in [0,1]$. The posterior density of $(\theta, \tau, \sigma^2)$ is proportional to
\begin{equation*}
\tilde{\ell} \big(\theta,\tau,\sigma^2 \, \big| \, \hat{\theta},\hat{\tau},\hat{\sigma}^2 \big)
\big(\xi \delta(\tau) + (1-\xi) \big) {\sigma}^{-2},
\end{equation*}
where $\tilde{\ell} \big(\theta,\tau,\sigma^2 \, \big| \, \hat{\theta},\hat{\tau},\hat{\sigma}^2 \big)$ is defined to be
\eqref{lik_prop}. Thus, the marginal posterior density of $\theta$ is proportional to
\begin{equation*}
\xi \, \int_0^{\infty} \int_{-\infty}^{\infty}
\tilde{\ell} \big(\theta, \tau, \sigma^2 \, | \, \hat{\theta},\hat{\tau},\hat{\sigma}^2 \big)
 \delta(\tau) \, {\sigma}^{-2} \, d \tau \, d \sigma^2 +
(1-\xi) \,  \int_0^{\infty} \int_{-\infty}^{\infty}
\tilde{\ell} \big(\theta, \tau, \sigma^2 \, | \, \hat{\theta},\hat{\tau},\hat{\sigma}^2 \big)
 \, {\sigma}^{-2} \, d \tau \, d \sigma^2.
\end{equation*}
It follows from the derivations presented in Appendices E and F that
this is equal to
\begin{equation}
\label{Mixture marginal posterior distribution}
\xi \, c_1(\hat{\tau}, \hat{\sigma}^2,2) \, f_{m+1}(\theta; \mu_1, {\sigma_1}^2(2))
+ (1-\xi)\, c_2(\hat{\sigma}^2,2) \, f_{m}(\theta; \hat{\theta}, \hat{\sigma}^2),
\end{equation}
where $\mu_1 = \hat{\theta}-\rho \hat{\tau}$, $\sigma_1^2(2) = (m\hat{\sigma}^2 + \hat{\tau}^2)(1-\rho^2)/(m+1)$ and
the functions $c_1$ and $c_2$ are defined by \eqref{c1} and \eqref{c_2}, respectively.
Therefore, the marginal posterior density of $\theta$ is
\begin{equation*}
\lambda(\hat{\sigma}, \hat{\tau}/\hat{\sigma}) \, f_{m+1}(\theta; \mu_1, {\sigma_1}^2(2))
 + \left(1-\lambda(\hat{\sigma}, \hat{\tau}/\hat{\sigma}) \right) \,f_{m}(\theta; \hat{\theta}, \hat{\sigma}^2),
\end{equation*}
where
$\lambda(\hat{\sigma}, \hat{\tau}/\hat{\sigma})
= \xi \, c_1(\hat{\tau}, \hat{\sigma}^2,2)/ \left(\xi \, c_1(\hat{\tau}, \hat{\sigma}^2,2)+(1-\xi)\, c_2(\hat{\sigma}^2,2) \right)$.
Thus
\begin{equation*}
\lambda(\hat{\sigma}, \hat{\tau}/\hat{\sigma})  =  \frac{1}{1 + k \, \hat{\sigma}\,{\left(m + \big(\hat{\tau}/\hat{\sigma}\big)^2 \right)}^{(m+1)/2}},
\end{equation*}
where
\begin{equation*}
k = \frac{(1-\xi)\,\sqrt{\pi}\,\,\Gamma(m/2)}{\xi\,m^{m/2} \, \Gamma((m+1)/2)}.
\end{equation*}

\bigskip

\noindent {\large{\bf {Appendix I: Properties of
$\boldsymbol{\big(\ell_B(\hat{\theta}, \hat{\tau}, \hat{\sigma}; \eta) - \hat{\theta} \big) / \hat{\sigma}}$ and \newline
$\boldsymbol{\big(u_B(\hat{\theta}, \hat{\tau}, \hat{\sigma}; \delta) - \hat{\theta} \big) / \hat{\sigma}}$
for the prior density \newline
$\boldsymbol{\pi(\theta, \tau, \sigma^2)  = \big(\xi \delta(\tau) + (1-\xi) \big) {\sigma}^{-2}}$}}}

\medskip

Suppose that the prior density is $\pi(\theta, \tau, \sigma^2)  = \big(\xi \delta(\tau) + (1-\xi) \big) {\sigma}^{-2}$
and that $\xi \in (0,1)$.
The marginal posterior density of $\theta$ is given by \eqref{first_marg_post_density_theta}. Therefore,
$\ell_B(\hat{\theta}, \hat{\tau}, \hat{\sigma}; \eta)$ is the solution for $v$ of
\begin{equation*}
\int_{-\infty}^v
\lambda(\hat{\sigma}, \hat{\tau}/\hat{\sigma})\, f_{m+1}(\theta; \mu_1, {\sigma_1}^2(2))
 + \left(1-\lambda(\hat{\sigma}, \hat{\tau}/\hat{\sigma}) \right) \, f_{m}(\theta; \hat{\theta}, \hat{\sigma}^2) \, d \theta
 = \eta
\end{equation*}
and $u_B(\hat{\theta}, \hat{\tau}, \hat{\sigma}; \delta)$ is the solution for $v$ of this equation, but with $\eta$ replaced
by $1 - \delta$. It follows from this that $\ell_B(\hat{\theta}, \hat{\tau}, \hat{\sigma}; \eta)$ is the solution for $v$ of
\begin{align*}
\lambda(\hat{\sigma}, \hat{\tau}/\hat{\sigma}) \,
&P \left( - \rho \frac{\hat{\tau}}{\hat{\sigma}} + \sqrt{\frac{(m + (\hat{\tau}/\hat{\sigma})^2)(1 - \rho^2)}{m+1}}
T_{m+1} \le \frac{v - \hat{\theta}}{\hat{\sigma}}\right) \\
&+ \left(1-\lambda(\hat{\sigma}, \hat{\tau}/\hat{\sigma}) \right) \,
P \left( T_m \le \frac{v - \hat{\theta}}{\hat{\sigma}}  \right)
 = \eta
\end{align*}
and $u_B(\hat{\theta}, \hat{\tau}, \hat{\sigma}; \delta)$ is the solution for $v$ of this equation, but with $\eta$ replaced
by $1 - \delta$, where $T_{m+1} \sim t_{m+1}$ and $T_{m} \sim t_{m}$. Thus
$\big(\ell_B(\hat{\theta}, \hat{\tau}, \hat{\sigma}; \beta) - \hat{\theta} \big) / \hat{\sigma}$ is the solution for $w$ of
\begin{align}
\label{eqn_appendix_I}
\lambda(\hat{\sigma}, \hat{\tau}/\hat{\sigma}) \,
&P \left( - \rho \frac{\hat{\tau}}{\hat{\sigma}} + \sqrt{\frac{(m + (\hat{\tau}/\hat{\sigma})^2)(1 - \rho^2)}{m+1}}
T_{m+1} \le w\right) \notag \\
&+ \left(1-\lambda(\hat{\sigma}, \hat{\tau}/\hat{\sigma}) \right) \,
P \left( T_m \le w \right)
 = \eta
\end{align}
and $\big(u_B(\hat{\theta}, \hat{\tau}, \hat{\sigma}; \delta) - \hat{\theta} \big) / \hat{\sigma}$
is the solution for $w$ of this equation, but with $\eta$ replaced
by $1 - \delta$. Clearly, \eqref{eqn_appendix_I} can be expressed in the following form
\begin{equation*}
\lambda(\hat{\sigma}, \hat{\tau}/\hat{\sigma})\,
F \left( \frac{w + \rho (\hat{\tau}/\hat{\sigma})}{\sqrt{(m + (\hat{\tau}/\hat{\sigma})^2) (1 - \rho^2)/(m+1)}} \, \Bigg | \,
t_{m+1} \right )
+ \left(1-\lambda(\hat{\sigma}, \hat{\tau}/\hat{\sigma}) \right) \,
F (w \, | \, t_m) = \eta,
\end{equation*}
where $F(\, \cdot \, | \, t_q)$ denotes the $t_q$ cumulative distribution function.
We see from this that both
$\big(\ell_B(\hat{\theta}, \hat{\tau}, \hat{\sigma}; \eta) - \hat{\theta} \big) / \hat{\sigma}$ and
$\big(u_B(\hat{\theta}, \hat{\tau}, \hat{\sigma}; \delta) - \hat{\theta} \big) / \hat{\sigma}$ are
functions of
$(\hat{\sigma}, \hat{\tau}/\hat{\sigma})$.

\bigskip

\noindent {\large{\bf {Appendix J: Plausibility argument leading to the prior \newline density
$\boldsymbol{\pi(\theta, \tau, \sigma^2)= \xi \delta(\tau)\sigma^{-2} + (1-\xi)\sigma^{-3}}$}}}

\medskip

Suppose that, conditional on $\sigma$, $\gamma = \tau / \sigma$ has the
improper prior density $\xi \delta(\gamma) + (1-\xi)$, where $\xi \in [0,1]$.
The prior probability that $\tau = 0$, conditional on $\sigma$, is the
same as the prior probability that $\tau / \sigma = 0$, conditional on $\sigma$. Thus the
prior probability that $\tau = 0$, conditional on $\sigma$, is $\xi$.
Now consider $0 < a < b$. The prior probability that $a \le \tau / \sigma \le b$, conditional on $\sigma$,
is $(b-a)(1-\xi)$. This implies that the prior probability that $a' \le \tau \le b'$, conditional on $\sigma$,
is $(b'-a')(1-\xi)/\sigma$, where $a' = a \sigma$ and $b' = b \sigma$.
A similar argument applies for $a < b < 0$.
Therefore, $\tau$ has prior density $\xi \delta(\tau) + (1-\xi)\sigma^{-1}$, conditional on $\sigma$.


\newpage

\noindent {\large{\bf {Appendix K: Marginal posterior distribution of
$\boldsymbol{\theta}$ for the prior density $\boldsymbol{\pi(\theta, \tau, \sigma^2)= \xi \delta(\tau)\sigma^{-g} + (1-\xi)\sigma^{-g-1}}$}}}

\medskip

Suppose that the prior density is $\pi(\theta, \tau, \sigma^2)  = \xi \delta(\tau){\sigma}^{-g} + (1-\xi) {\sigma}^{-g-1}$,
where $m+g > 0$.
Also suppose that $\xi \in [0,1]$. The posterior density of $(\theta, \tau, \sigma^2)$ is proportional to
\begin{equation*}
\tilde{\ell} \big(\theta,\tau,\sigma^2 \, \big| \, \hat{\theta},\hat{\tau},\hat{\sigma}^2 \big)
\, \big(\xi \delta(\tau){\sigma}^{-g} + (1-\xi) {\sigma}^{-g-1} \big),
\end{equation*}
where $\tilde{\ell} \big(\theta,\tau,\sigma^2 \, \big| \, \hat{\theta},\hat{\tau},\hat{\sigma}^2 \big)$ is defined to be
\eqref{lik_prop}. Thus, the marginal posterior density of $\theta$ is proportional to
\begin{equation*}
\xi \, \int_0^{\infty} \int_{-\infty}^{\infty}
\tilde{\ell} \big(\theta, \tau, \sigma^2 \, | \, \hat{\theta},\hat{\tau},\hat{\sigma}^2 \big)
 \delta(\tau) \, {\sigma}^{-g} \, d \tau \, d \sigma^2 +
(1-\xi) \,  \int_0^{\infty} \int_{-\infty}^{\infty}
\tilde{\ell} \big(\theta, \tau, \sigma^2 \, | \, \hat{\theta},\hat{\tau},\hat{\sigma}^2 \big)
 \, {\sigma}^{-g-1} \, d \tau \, d \sigma^2.
\end{equation*}
It follows from the derivations presented in Appendices E and F that
this is equal to
\begin{equation*}
\xi \, c_1(\hat{\tau}, \hat{\sigma}^2, g) \, f_{m+g-1}(\theta; \mu_1, {\sigma_1}^2(g))
+ (1-\xi)\, c_2(\hat{\sigma}^2, g+1) \, f_{m+g-1}(\theta; \hat{\theta}, {\sigma_2}^2(g+1)),
\end{equation*}
where $\mu_1 = \hat{\theta}-\rho \hat{\tau}$, $\sigma_1^2(g) = (m\hat{\sigma}^2 + \hat{\tau}^2)(1-\rho^2)/(m+g-1)$,
$c_1(\hat{\tau}, \hat{\sigma}^2, g)$ is given by \eqref{c1} and, in accordance with the definitions at the end of Appendix F,
\begin{equation*}
c_2(\hat{\sigma}^2, g+1) = 2\pi\sqrt{1-\rho^2} {\left(\frac{m \hat{\sigma}^2}{2}\right)}^{-(m+g-1)/2} \Gamma\left(\frac{m+g-1}{2}\right),
\end{equation*}
and
\begin{equation*}
{\sigma_2}^2(g+1) = \frac{m}{m+g-1} \, \hat{\sigma}^2.
\end{equation*}
Therefore, the marginal posterior density of $\theta$ is
\begin{equation}
\label{marginal_post_density_2nd_general}
\tilde{\lambda}(\hat{\tau}/\hat{\sigma},g) \, f_{m+g-1}(\theta; \mu_1, {\sigma_1}^2(g))
 + \left(1-\tilde{\lambda}(\hat{\tau}/\hat{\sigma},g) \right) \,f_{m+g-1}(\theta; \hat{\theta}, {\sigma_2}^2(g+1)),
\end{equation}
where
$\tilde{\lambda}(\hat{\tau}/\hat{\sigma},g)
= \xi \, c_1(\hat{\tau}, \hat{\sigma}^2, g)/ \left(\xi \, c_1(\hat{\tau}, \hat{\sigma}^2, g)+(1-\xi)\, c_2(\hat{\sigma}^2, g+1) \right)$.
Thus
\begin{equation*}
\tilde{\lambda}(\hat{\tau}/\hat{\sigma},g)
=  \frac{1}{1 + \tilde{k}(g) \, {\left(m + \big(\hat{\tau}/\hat{\sigma}\big)^2 \right)}^{(m+g-1)/2}},
\end{equation*}
where
\begin{equation*}
\tilde{k}(g) = \frac{\sqrt{2 \pi}\,(1-\xi)}{\xi \, m^{(m+g-1)/2}}.
\end{equation*}

\newpage

\noindent {\large{\bf {Appendix L: Properties of
$\boldsymbol{\big(\ell_B(\hat{\theta}, \hat{\tau}, \hat{\sigma}; \eta) - \hat{\theta} \big) / \hat{\sigma}}$ and \newline
$\boldsymbol{\big(u_B(\hat{\theta}, \hat{\tau}, \hat{\sigma}; \delta) - \hat{\theta} \big) / \hat{\sigma}}$
for the prior density \newline
$\boldsymbol{\pi(\theta, \tau, \sigma^2)= \xi \delta(\tau)\sigma^{-g} + (1-\xi)\sigma^{-g-1}}$}}}

\medskip

Suppose that the prior density is $\pi(\theta, \tau, \sigma^2)  = \xi \delta(\tau){\sigma}^{-g} + (1-\xi) {\sigma}^{-g-1}$,
where $m+g > 0$.
Also suppose that $\xi \in [0,1]$.
The marginal posterior density of $\theta$ is given by \eqref{marginal_post_density_2nd_general}. Therefore,
$\ell_B(\hat{\theta}, \hat{\tau}, \hat{\sigma}; \eta)$ is the solution for $v$ of
\begin{equation*}
\int_{-\infty}^v
\tilde{\lambda}(\hat{\tau}/\hat{\sigma},g) \, f_{m+g-1}(\theta; \mu_1, {\sigma_1}^2(g))
 + \left(1-\tilde{\lambda}(\hat{\tau}/\hat{\sigma},g) \right) \,f_{m+g-1}(\theta; \hat{\theta}, {\sigma_2}^2(g+1)) \, d \theta
 = \eta
\end{equation*}
and $u_B(\hat{\theta}, \hat{\tau}, \hat{\sigma}; \delta)$ is the solution for $v$ of this equation, but with $\eta$ replaced
by $1 - \delta$. It follows from this that $\ell_B(\hat{\theta}, \hat{\tau}, \hat{\sigma}; \beta)$ is the solution for $v$ of
\begin{align*}
\tilde{\lambda}(\hat{\tau}/\hat{\sigma},g) \,
&P \left( - \rho \frac{\hat{\tau}}{\hat{\sigma}} + \sqrt{\frac{(m + (\hat{\tau}/\hat{\sigma})^2)(1 - \rho^2)}{m+g-1}} \,
T_{m+g-1} \le \frac{v - \hat{\theta}}{\hat{\sigma}}\right) \\
&+ \left(1-\tilde{\lambda}(\hat{\tau}/\hat{\sigma},g)\right) \,
P \left( \sqrt{\frac{m}{m+g-1}} \, T_{m+g-1} \le \frac{v - \hat{\theta}}{\hat{\sigma}}  \right)
 = \eta
\end{align*}
and $u_B(\hat{\theta}, \hat{\tau}, \hat{\sigma}; \delta)$ is the solution for $v$ of this equation, but with $\eta$ replaced
by $1 - \delta$, where $T_{m+g-1} \sim t_{m+g-1}$. Thus
$\big(\ell_B(\hat{\theta}, \hat{\tau}, \hat{\sigma}; \eta) - \hat{\theta} \big) / \hat{\sigma}$ is the solution for $w$ of
\begin{align*}
\label{eqn_appendix_I}
\tilde{\lambda}(\hat{\tau}/\hat{\sigma},g) \,
&P \left( - \rho \frac{\hat{\tau}}{\hat{\sigma}} + \sqrt{\frac{(m + (\hat{\tau}/\hat{\sigma})^2)(1 - \rho^2)}{m+g-1}} \,
T_{m+g-1} \le w \right) \\
&+ \left(1-\tilde{\lambda}(\hat{\tau}/\hat{\sigma},g)\right) \,
P \left( \sqrt{\frac{m}{m+g-1}} \, T_{m+g-1} \le w \right)
 = \eta
\end{align*}
and $\big(u_B(\hat{\theta}, \hat{\tau}, \hat{\sigma}; \delta) - \hat{\theta} \big) / \hat{\sigma}$
is the solution for $w$ of this equation, but with $\eta$ replaced
by $1 - \delta$. Clearly, this equality can be expressed in the following form
\begin{align*}
\tilde{\lambda}(\hat{\tau}/\hat{\sigma},g)\,
&F \left( \frac{w + \rho (\hat{\tau}/\hat{\sigma})}{\sqrt{(m + (\hat{\tau}/\hat{\sigma})^2) (1 - \rho^2)/(m+g-1)}} \, \Bigg | \,
t_{m+g-1} \right ) \\
&+ \left(1-\tilde{\lambda}(\hat{\tau}/\hat{\sigma},g) \right) \,
F \left(\sqrt{\frac{m+g-1}{m}} \, w \, \Bigg| \, t_{m+g-1} \right) = \eta,
\end{align*}
where $F(\, \cdot \, | \, t_q)$ denotes the $t_q$ cumulative distribution function.
We see from this that both
$\big(\ell_B(\hat{\theta}, \hat{\tau}, \hat{\sigma}; \eta) - \hat{\theta} \big) / \hat{\sigma}$ and
$\big(u_B(\hat{\theta}, \hat{\tau}, \hat{\sigma}; \delta) - \hat{\theta} \big) / \hat{\sigma}$ are
functions of
$\hat{\tau}/\hat{\sigma}$.



\newpage

\noindent {\bf References}

\smallskip


\rf Box, G.E.P., Tiao, G.C., 1973. Bayesian Inference in Statistical Analysis. Wiley, New York.

\smallskip

\rf Chipman, H., George, E.I., McCulloch, R.E. (2001). The Practical Implementation of Bayesian Model Selection (with discussion). IMS Lecture Notes - Monograph Series, 38, 65 - 116.

\smallskip

\rf Farchione, D., Kabaila, P., 2008. Confidence intervals for the normal mean utilizing uncertain prior information.
Statistics and Probability Letters 78, 1094--1100.

\smallskip

\rf Ferentinos, K.K., Karakostas, K.X., 2006. More on shortest and equi tails confidence intervals. Communications in Statistics
-Theory and Methods 35, 821--829.

\smallskip

\rf Johnstone, I.M., Silverman, B.W., 2004. Needles and straw in haystacks: empirical Bayes estimates of possibly sparse sequences.
Annals of Statistics 32, 1594--1649.

\smallskip

\rf Johnstone, I.M., Silverman, B.W., 2005. Bayes selection of wavelet thresholds.
Annals of Statistics 33, 1700--1752.

\smallskip

\rf Kabaila, P., 2009. The coverage properties of confidence regions after model selection. International Statistical Review, 77,
405--414.

\smallskip

\rf Kabaila, P., Giri, K., 2009. Confidence intervals in regression utilizing uncertain prior information.
Journal of Statistical Planning and Inference 139, 3419--3429.

\smallskip

\rf Kabaila, P., Giri, K., Leeb, H., 2010. Admissibility of the usual confidence interval in linear regression. Electronic Journal of Statistics
4, 300--312.

\smallskip

\rf Kabaila, P., Giri, K., 2013. Further properties of frequentist confidence intervals in regression that utilize uncertain prior information.
Australian \& New Zealand Journal of Statistics, 55, 259--270.

\smallskip

\rf Miller, A., 2002. Subset Selection in Regression, Second Edition. Chapman \& Hall/CRC.

\smallskip

\rf Mitchell, T.J., Beauchamp, J.J., 1988. Bayesian variable selection in linear regression.
Journal of the American Statistical Association, 83, 1023--1032.

\smallskip

\rf O'Hara, R.B., Sillanp\"a\"a, M.J., 2009. A review of Bayesian variable selection methods: what, how and which. Bayesian Analysis
4, 85--118.

\smallskip

\rf Zellner, A., 1986. On assessing prior distributions and Bayesian regression analysis with $g$-prior distributions.
Pages 233--243 of Bayesian Inference and Decision Techniques, Essays in Honor of Bruno de Finetti, edited by Prem K. Goel
and Arnold Zellner, Elsevier, Amsterdam, The Netherlands.



\end{document}